\documentclass[12pt]{article}

\usepackage{amsfonts}
\usepackage{amsmath}
\usepackage{amssymb}
\usepackage{graphicx}

\def\dubluunu{\mbox{$1\hspace*{-0,7mm}\rule{0,2mm}{2,7mm}\,$}}
\def\PP{{\mathbb{P}}}
\def\E{{\mathbb{E}}}
\def\rr{{\mathbb{R}}}
\def\nn{{\mathbb{N}}}
\def\mk{\medskip }
\def\bk{\bigskip }
\def\sk{\smallskip }
\def\ct{con\-ti\-nuous}
\def\vsp{\vspace*{1,5mm}\\ }
\def\q{\quad }
\def\qq{\qquad }

\def\sumk{\dd\sum^m_{k=1}}
\def\oo{{\omega}}
\def\ooo{{\Omega}}
\def\<{\left<}
\def\>{\right>}
\def\({\left(}
\def\){\right)}
\def\ff{\forall }
\def\D{\Delta }
\def\9{{\infty}}
\def\barr{\begin{array}}
\def\earr{\end{array}}
\def\ov{\overline}
\def\wt{\widetilde}
\def\vp{{\varepsilon}}

\def\pp{{\partial}}
\def\dd{\displaystyle}
\def\a{{\alpha}}
\def\b{{\beta}}
\def\vf{{\varphi}}

\def\pas{\mathbb{P}\mbox{-a.s.}}
\def\de{{\delta}}

\def\cald{{\mathcal{D}}}
\def\calf{{\mathcal{F}}}

\def\calo{\mathcal{O}}
\def\calf{\mathcal{F}}

\def\na{{\nabla}}
\def\ifr{\mbox{ if }}
\def\inr{\mbox{ in }}
\def\onr{\mbox{ on }}
\def\fwg{following}
\def\eq{equation}

\def\beq#1{\begin{equation}\label{e#1}}
\def\eeq{\end{equation}}
\def\1{^{-1}}

\newtheorem{theorem}{Theorem}[section]

\newtheorem{definition}[theorem]{Definition}

\title{\bf The stochastic reflection problem with multiplicative noise\footnote{This work is supported by CNCSIS project PN ID-/2011.}}
\author{Viorel Barbu}
\date{Al.I. Cuza University and Octav Mayer Institute\\ of the Romanian Academy, Ia\c si, Romania}

\begin{document}

\maketitle

\begin{abstract}This paper addresses the existence and uniqueness of strong solutions to the stochastic variational inequality
$$dX-\D X\,dt+F(t,\xi,X)dt+\b(X)dt\ni\dd\sum^m_{k=1}X\mu_k d\b_k(t)+f(t)dt$$ in a bounded domain $\calo\subset \rr^d$ with Dirichlet homogeneous conditions. Here $\b(r)=0$ for $r>0$, $\b(0)=]-\9,0]$, $\b(r)\ne0$ for $r<0$. An application to the one-phase Stefan problem with a stochastic heat source is given. One studies also the corresponding stochastic parabolic equation with Signorini boundary conditions $\dd\frac{\pp X}{\pp\nu}+\b(X)\ni0$ on $\pp\calo$.\bk
\end{abstract}

\section{Introduction}

This work is concerned with the existence and uniqueness of the stochastic variational inequality

\begin{equation}\label{e1.1}
\left\{\barr{l}
d_t X(t,\xi)-\D_\xi X(t,\xi)dt+F(t,\xi,X(t,\xi))dt\vsp
\qquad+\b(X(t,\xi))dt\ni\dd\sum^m_{k=1}X(t,\xi)\mu_k(t,\xi)d\b_k(t)+f(t,\xi)dt\\\hfill \mbox{in }(0,T)\times\calo=Q_T,\vsp
X(0,\xi)=x(\xi),\ \xi\in\calo,\vsp
X(t,\xi)=0\ \mbox{ on }(0,T)\times\pp\calo=\Sigma_T,\earr\right.\end{equation}
where $\calo$ is an open and bounded domain of $\rr^d$, $d\ge1$, with smooth boundary $\pp\calo$, $\{\b_k(t)\}^m_{k=1}$ is a system of mutually independent Brownian motions on a probability space $\{\ooo,\calf,\PP\}$ and $\b:\rr\to2^\rr$ is the multivalued function
\begin{equation}\label{e1.2}
\b(r)=0\ \mbox{ for }r>0,\qquad\b(0)=]-\9,0],\qq \b(r)=\emptyset\ \mbox{ for }r<0.\end{equation}
As regard the functions $F:Q_T\times\rr\to\rr$, $\mu_k:[0,T]\times\calo\to\rr$, we assume that the \fwg\ hypotheses hold.\mk

(i) $F=F(t,\xi,r)$ is measurable in $(t,\xi)\in Q_T$, $F(t,\xi,0)=0$ and
\begin{equation}\label{e1.3}
|F(t,\xi,r)-F(t,\xi,\bar r)|\le\a|r-\bar r|,\ \ff r,\bar r\in\rr,\ \mbox{a.e.}\ (t,\xi)\in Q_T\end{equation}
where $\a>0.$\sk

(ii) $\mu_k\in C^2([0,T]\times\ov\calo),\ k=1,...,m.$\mk

The function $f:Q_T\times\ooo\to\rr$ is an adapted process from $(0,T)$ to $L^2(\calo)$ with respect to the probability basis $\{\ooo,\calf,\{\calf_t\},\PP\},$ which will be made precise later on.

Formally, problem \eqref{e1.1} can be rewritten as a free boundary value problem (the stochastic reflection problem)
\begin{equation}\label{e1.5}
\left\{\!\!\!\barr{l}
d_tX-\D_\xi X\,dt+F(t,\xi,X)dt=f\,dt+\dd\sum^m_{k=1}X\mu_kd\b_k\\\hfill\mbox{ in }[(t,\xi)\in Q_T;\,X(t,\xi)>0]\vsp
X\ge0,\ d_t X-\D_\xi X\,dt+F(t,\xi,X)dt\ge f\,dt+\dd\sum^m_{k=1}X\mu_kd\b_k\mbox{ in }Q_T\vsp
X(0)=x\mbox{ in }\calo,\ \ \ X=0\mbox{ on }\Sigma_T.\earr\right.\end{equation}
As regards the literature on stochastic parabolic variational inequalities of this form or, more generally, on parabolic differential equations driven by Gaussian noise and reflected to the boundary of a convex set of $\rr^d$, the works \cite{6}, \cite{8}, \cite{10}, \cite{11}, \cite{13} must be cited. (See also \cite{2}, \cite{3}, \cite{7} for the finite-dimensional case and \cite{4}, \cite{5}, \cite{14} or \cite{1}, Chapter 1, for existence theory in infinite-dimensional Hilbert spaces.) It should be said, however, that all refer to the case of stochastic equations with additive cylindrical Gaussian noise and the solution $X$ is taken in a generalized sense. Roughly speaking, such a solution $X$ satisfies \eqref{e1.1} where $\b(X)$ is replaced by a measure on $(0,T)\times\calo$ with the support in $[X=0]$. Theorem \ref{t2.2} below provides the existence and uniqueness of a strong solution to \eqref{e1.1} which is more appropriate of classical formulation \eqref{e1.5} of the obstacle parabolic variational inequality. We develop here a parallel study for the stochastic parabolic equation \eqref{e1.1} with Signorini boundary value conditions. Namely,
\begin{equation}\label{e1.6}
\left\{\barr{l}
d_tX(t,\xi)-\D_\xi X(t,\xi)dt+F(t,\xi,X(t,\xi))dt\vsp
\qq\qq=f(t,\xi)dt+\dd\sum^m_{k=1}X(t,\xi)\mu_kd\b_k(t)\ \mbox{ in }Q_T\vsp
X(0,\xi)=x(\xi),\ \ \xi\in\calo\vsp
\dd\frac{\pp  X}{\pp\nu}\ (t,\xi)+\b(X(t,\xi))\ni0\ \mbox{ on }\Sigma_T,\earr\right.\end{equation}
where $\dd\frac{\pp}{\pp\nu}$ is the outward normal derivative to $\pp\calo$.

As made precise in Section 4, problem \eqref{e1.1} can be used to describe the dynamic of the one-phase Stefan melting (solidification) problem in presence of a stochastic heat source of the form $\dd\sum^m_{k=1}\theta\mu_k d\b_k$, where $\theta =\theta(t,\xi)$ is the temperature. As regards \eqref{e1.6}, it is a stochastic parabolic equations with unilateral conditions on the boundary.

\section{The main results}
\setcounter{equation}{0}

\subsection*{Notation}

Given the stochastic basis $\{\ooo,\PP,\calf,\{\calf_t\}_{t\ge0}\}$ and $Z$ a Hilbert space with the norm $\|\cdot\|_Z,$ we denote by $M^p_\PP(0,T;Z)$ the space of all progressively measurable $Z$-valued processes $X:\ooo\times(0,T)\to Z$ with respect to the filtration $\{\calf_t\}_{t\ge0}$   such that
$$\|X\|^p_{M^p_\PP(0,T;Z)}=\E\int^T_0\|X(t)\|^p_Zdt<\9,$$where $\E$ is the expectation. Denote by $C_\PP([0,T];Z)$ the space of all processes $X\in  C([0,T];L^2(\ooo,\PP,\calf,Z))$. By $L^p_{ad}(\ooo;C([0,T];L^2(\calo))$ we denote the space of all adapted processes with respect to $\{\calf_t\}$, $X:\ooo\to C([0,T];L^2(\calo))$ such that
$$\E\|X\|^p_{C([0,T];L^2(\calo))}<\9.$$By $L^p(\calo)$, where $\calo$ is a bounded and open subset of $\rr^d$ and $1\le p\le\9$, we denote the space of all $p$-summable functions on $\calo$ with the norm denoted by $|\cdot|_p$. By $L^p(0,T;Z)$ we denote the space of all Bochner $Z$-valued $p$-summable functions on $(0,T)$ and by $W^{1,p}([0,T];Z)$ the space of all absolutely \ct\ functions $y:[0,T]\to Z$ such that $\dd\frac{dy}{dt}\in L^p(0,T;Z)$. The norm of $Z$ is denoted by $\|\cdot\|_Z$.

In the \fwg, $H^1(\calo)$, $H^1_0(\calo)$ and $H^2(\calo)$ are standard Sobolev spaces on $\calo$. We also denote by $H^{-1}(\calo)$ the dual space of $H^1_0(\calo)$ and set
$$H^{2,1}(Q_T)=\left\{y\in L^2(0,T;H^1_0(\calo)\cap H^2(\calo),\frac{\pp y}{\pp t}\in L^2(Q_T)\right\}.$$

\begin{definition}\label{d2.1} {\rm By strong solution $X$ to \eqref{e1.1}, we mean a process $X\in C_\PP([0,T];L^2(\calo))\cap M^2_\PP(0,T;H^1_0(\calo))$ such that   there is $\eta\in M^2_\PP(0,T;L^2(\calo))$ and $\pas$,
\begin{equation}\label{e2.1}
\left\{\barr{l}
\eta(t,\xi)\in\b(X(t,\xi))\q\mbox{a.e. }(t,\xi)\in Q_T\vsp
X(t)=x+\!\!\dd\int^t_0\!\D_\xi X(x)ds-\!\!\int^t_0\!F(s,\xi,X(s))+\eta(s))ds\vsp
\qquad+\!\!\dd\int^t_0\!f(s)ds+\!\dd\sum^N_{k=1}\int^t_0\!\mu_k(s)X(s)d\b_k(s),\ t\in [0,T],\earr\right.\end{equation}
where the integral is taken in the sense of Ito in the space $H^{-1}(\calo)$. (We refer to \cite{9} for existence and uniqueness of such a solution.)
}\end{definition}
It should be noted that, for any such a solution $X$ to \eqref{e1.1} in the sense of Definition \ref{d2.1}, the Ito's formula is applicable.

Taking into account that the relation $\eta\in\b(X)$ is equivalent with
\begin{equation}\label{e2.2}
X\ge0,\ \ \eta\ge0,\ \  X\eta=0,\ \ \mbox{ a.e. in }Q_T,\end{equation}
it follows by \eqref{e2.1} that $X$ defined above is a solution to \eqref{e1.1} in the sense of~\eqref{e1.2}.

\begin{theorem}\label{t2.2} Let $f\in M^2_\PP(0,T);L^2(\calo))$ and $x\in H^1_0(\calo)$ be such that $x\ge0$, a.e. in $\calo$. Then \eqref{e1.1} has a unique strong solution $X$. Moreover, we have
\begin{equation}\label{e2.7a} X\in L^p_{ad}(\ooo;C([0,T];L^2(\calo))\cap M^p_\PP(0,T;H^2(\calo)),\end{equation}for all $1\le p<2.$
\end{theorem}

In particular, it follows by \eqref{e2.7a} that $\D_\xi X\in L^2(Q_T)$, $\pas$, and so  \eqref{e2.1} holds, a.e. on $(0,T)\times\calo\times\ooo$.

Now, we come back to \eqref{e1.6} to define a strong solution for this equation. To this purpose, we define the multivalued nonlinear operator $A(t):V\to V'$, $V=H^1(\calo)$,
\begin{equation}\label{e2.3}
\barr{lcl}
_{V'}\<A(t)X,\vf\>_V&=&\dd\int_\calo(\na X(\xi)\cdot\na\vf(\xi)+F(t,\xi,X(\xi))\vf(\xi))d\xi\vsp
&&+\dd\int_{\pp\calo}\eta(X)(\xi)\vf(\xi)d\xi,\ \ff\vf\in V,\ t\in[0,T],\earr\end{equation}
where $\eta\in L^2(\pp\calo)$ and
\begin{equation}\label{e2.4}
\eta(X)(\xi)\in\b(X(\xi)),\ \ \mbox{a.e. }\xi\in\pp\calo.\end{equation}
Here $V'$ is the dual of $V$ and $_{V'}\<\cdot,\cdot\>_V$ is the duality between $V$ and $V'$ induced by the space $H=L^2(\calo)$. (We have $V\subset H\subset V'$ in algebraic and topological sense.)

The operator $A$ is not everywhere defined on $V$ but on the set $D(A)$ of all $X\in V$ for which there is $\eta\in L^2(\pp\calo)$ satisfying \eqref{e2.4}. Roughly speaking, $A(t)=-\D$ with the nonlinear boundary conditions $\dd\frac{\pp X}{\pp\nu}+\b(X)\ni0$ on $\pp\calo$.

\begin{definition}\label{d2.3} {\rm $X$ is said to be strong solution to \eqref{e1.6} if the \fwg\ conditions hold:
\begin{equation}\label{e2.5}
X\in C_\PP([0,T];L^2(\calo))\cap M^2_\PP(0,T;V),\qq\qq\qq\q\end{equation}
\begin{equation}\label{e2.6}
\barr{lcl}
\dd X(t)-x\in &-&\dd\int^t_0A(s)X(s)ds+\dd\int^t_0f(s)ds\vsp
&+&\dd\int^t_0X(s)\dd\sum^m_{k=1}\mu_k(s)d\b_k(s),\  \ff t\in[0,T],\ \pas\earr\end{equation}
}\end{definition}
Here $\dd\int^t_0 A(s)X(s)ds$ is the Bochner $V'$-valued integral of the multivalued function $t\to A(t) X(t)$ and $\dd\int^t_0X(s)\dd\sum^m_{k=1}\mu_k(s)d\b_k(s)\in L^2(\calo)$, $\pas$, is taken in the Ito sense.

We may rewrite \eqref{e2.6} as the multivalued stochastic integral equation in~$V'$
\begin{equation}\label{e2.7}
\barr{l}
dX(t)+A(t)X(t)dt\ni f(t)+\dd\sum^m_{k=1}X(t)\mu_k d\b_k(t),\ t\in(0,T),\vsp
X(0)=x.\earr\end{equation}

\begin{theorem}\label{t2.4} Let $f\in M^2_\PP(0,T);L^2(\calo))$ and $x\in H^1_0(\calo)$ be such that $x\ge0$, a.e. in $\calo$. Then \eqref{e1.6} has a unique strong solution $X$. Moreover, $X\ge0$ a.e. in $(0,T)\times\calo\times\ooo$.
\end{theorem}

In particular, it follows by \eqref{e2.7a} that \eqref{e2.6} holds, a.e. on $(0,T)\times\calo\times\ooo$.

\section{Proofs}
\setcounter{equation}{0}

\subsection{Proof of Theorem \ref{t2.2}}

We set
$$\mu(t,\xi)=\sumk\mu_k(t,\xi)\b_k(t),\q(t,\xi)\in Q_T$$and$$\wt\mu(t,\xi)=\sumk\(\frac{\pp\mu_k}{\pp t}\,(t,\xi)\b_k(t)+\frac12\,\mu^2_k(t,\xi)\),\q(t,\xi)\in Q_T.$$By the substitution
\begin{equation}\label{e3.1}
X(t,\xi)=e^{\mu(t,\xi)}y(t,\xi),\ \oo\in\ooo,\end{equation}\eq\ \eqref{e1.1} reduces to the random parabolic \eq
\begin{equation}\label{e3.2}
\left\{\barr{l}
\dd\frac{\pp y}{\pp t}\,(t,\xi)-e^{-\mu(t,\xi)}\D(e^{\mu(t,\xi)}y(t,\xi))+e^{\mu(t,\xi)}F(t,\xi,e^{\mu(t,\xi)}y(t,\xi))\vsp
\qq\q+\wt\mu(t,\xi)y(t,\xi)+\b(y(t,\xi))\ni f(t,\xi),\q(t,\xi)\in Q_T,\vsp
y(0,\xi)=x(\xi)\q\inr\ \calo,\vsp
y(t,\xi)=0\q\onr\ \Sigma_T.\earr\right.\end{equation}Indeed, if $y=y(t,\xi)$ is a sufficiently regular solution to equation \eqref{e3.2}, then, by Ito's formula, we have$$dX=yd(e^\mu)+e^\mu\ \dd\frac{\pp y}{\pp t}\,dt\q\inr\ Q_T$$and$$\barr{lcl}d(e^\mu)&=&e^\mu d\mu+\dd\frac12\,e^\mu\dd\sum^N_{k=1}\mu^2_kdt\vsp
&=&e^\mu\dd\sum^N_{k=1}\mu_kd\b_k+\dd\frac12\,e^\mu\dd\sum^N_{k=1}\mu^2_kdt+e^\mu\dd\sum^N_{k=1}\dd\frac{\pp\mu_k}{\pp t}\,\b_kdt.\earr$$Then, substituting $X$ into \eqref{e1.1}, we obtain \eqref{e3.2}, as claimed.\bk

Now, by \eqref{e3.2}, we have $\pas$
\begin{equation}\label{e3.3}
\left\{\barr{l}
\dd\frac{\pp y}{\pp t}-\D y+\wt F(t,y)+g\cdot\na y+\b(y)\ni f\q\inr\ Q_T,\vsp
y(0)=x\q\inr\ \calo,\vsp
y=0\q\onr\ \Sigma_T,\earr\right.\end{equation}where
\begin{equation}\label{e3.4}
\left\{\barr{l}
\wt F(t,y)=e^{-\mu(t,\xi)}F(t,\xi,e^{\mu(t,\xi)}y)+\wt\mu y-(|\na\mu|^2+\D\mu)y,\vsp
g=-2\na\mu.\earr\right.\end{equation}

Any regular and predictible  solution $t\to y(t)$ with respect to the stochastic pro\-ba\-bi\-lity basis $\{\ooo,\calf,\calf_t,\PP\}$ to \eq\ \eqref{e3.3} leads via transformation \eqref{e3.1} to a solution $X$ to \eqref{e1.1} in the sense of the above definition (see, e.g., \cite{9}, p.~232).

We note that \eqref{e3.3} is a random parabolic variational inequality for which the standard existence theory applies (see, e.g., \cite{1}, p.~209). Namely, for each $\oo\in\ooo$, \eqref{e3.3} has a unique solution $y=y(t,\xi,\oo),$
\begin{equation}\label{e3.5}
y\in W^{1,2}([0,T];L^2(\calo))\cap L^2(0,T;H^1_0(\calo)\cap H^2(\calo)),\end{equation}which is obtained as limit for $\vp\to0$ of the solution $y_\vp$ to the penalized \eq
\begin{equation}\label{e3.6}
\barr{l}
\dd\frac{\pp y}{\pp t}-\D y_\vp+\wt F(t,y_\vp)+g\cdot\na y_\vp+\b_\vp(y_\vp)=f\q\inr\ Q_T,\vsp
y_\vp(0)=x\q\inr\ \calo,\qq y_\vp=0\q\onr\ \Sigma_T,\earr\end{equation}where $\b_\vp(r)=\dd\frac1\vp\,(r-(1+\vp\b)^{-1}r)=-\dd\frac1\vp\,r^-,$ $\ff r\in\rr.$

Indeed, if $y_\vp\in W^{1,2}([0,T];L^2(\calo))\cap L^2(0,T;H^1_0(\calo)\cap H^2(\calo))$ is the solution to \eqref{e3.6} taking into account that, by \eqref{e1.3}, a.e., $(t,\xi)\in Q_T$,
\begin{equation}\label{e3.7}
\barr{l}
|\wt F(t,r)|\le\a|r|,\q\ff r\in\rr,\ (t,\xi)\in Q_T,\vsp
|g|\le C\sup\{|\b_k(t)|,\ t\in[0,T],\ k=1,...,m\}=\de(\oo),\ \oo\in\ooo,\earr\end{equation}after some straightforward calculation, we obtain the estimate
\begin{equation}\label{e3.8}
|y_\vp(t)|^2_2+\dd\int^t_0\|y_\vp(s)\|^2_{H^1_0(\calo)}ds\le|x|^2_2+\dd\int^t_0|f(s)|^2_2ds+\de^2(\oo),\ t\in[0,T].\end{equation}

Similarly, multiplying \eqref{e3.6} by $\b_\vp(y_\vp)$ and integrating on $(0,t)\times\calo$, we obtain
\begin{equation}\label{e3.9}
\barr{r}\dd\int^t_0(|\b_\vp(y_\vp(s))|^2_2+|\D y_\vp(s)|^2_2)ds\le C\(\dd\int^t_0|f(s)|^2_2ds+t\de^2(\oo)\),\vsp \ff t\in[0,T],\earr\end{equation}because
$$\dd\int_\calo\dd\frac{\pp y_\vp}{\pp t}\,(t)\b_\vp(y_\vp(t))d\xi=\dd\frac1{2\vp}\ \frac d{dt}\ |y^-_0(t)|^2_2,\mbox{ a.e. }t\in[0,T],\mbox{ and }x^-=0\ \inr\ \calo_0.$$(Everywhere in the \fwg, we denote by $C$ several constants independent of $\vp$ and $\oo\in\ooo.$) Taking into account that $\E(\de^2)<\9$, we obtain by \eqref{e3.8}, \eqref{e3.9} that

\begin{equation}\label{e3.10}
\barr{l}
\E\left[|y_\vp(t)|^2_2\ +\dd\int^T_0(\|y_\vp(t)\|^2_{H^1_0(\calo)}
\right.\vsp
\qquad\qquad\quad\left.+\left\|\dd\frac{\pp y_\vp}{\pp t}\,(t)\right\|_2+\|\D y_\vp(t)\|^2_2+|\b_\vp(y_\vp(t))|^2_2))dt\right]\vsp
\qquad\qquad\quad\le C\(|x|^2_2+\E\dd\int^T_0|f(t)|^2_2dt\),\ \ff t\in[0,T].\earr\end{equation}
We note also that, since $y_\vp$ is obtained by an iteration process given by the contraction principle, $y_\vp(t,\oo)$ is progressively measurable and, therefore, adapted to the filtration $\{\calf_t\}_{t\ge0}$.

Moreover, by \eqref{e3.6} we have also for $\vp,\vp'>0$

$$\barr{l}
\dd\frac12\,|y_\vp(t)-y_{\vp'}(t)|^2_2
+\dd\int^t_0|\na(y_\vp(s)-y_{\vp'}(s)|^2ds\vsp\qq\qq+
\dd\int^t_0\int_\calo(\b_\vp(y_\vp(s))-
\b_{\vp'}(y_\vp(s)))y_\vp(s)-y_{\vp'}(s))ds\,d\xi\vsp
\qq\qq\le\dd\frac12\int^t_0\int_\calo|{\rm div}\ g|\,|y_\vp(s)-y_{\vp'}(s))|^2ds\,d\xi\vsp
\qq\qq\le\sup\{|{\rm div}\ g|;\ (t,\xi)\in Q\}\dd\int^t_0|y_\vp(s)-y_{\vp'}(s)|^2_2ds\earr$$
and
$$
(\b_\vp(y_\vp)-\b_{\vp'}(y_\vp))(y_\vp-y_{\vp'})
\ge(\vp\b_\vp(y_\vp)-\vp'\b_{\vp'}(y_{\vp'}))
(\b_\vp(y_\vp)-\b_{\vp'}(y_{\vp'})),$$and so, by \eqref{e3.8}, \eqref{e3.9}, we have
\begin{equation}\label{e3.11}
|y_\vp(t)-y_{\vp'}(t)|^2_2+\dd\int^t_0|\na(y_\vp(s)-y_{\vp'}(s)|^2_2ds\le C_1(\oo)(\vp+\vp'),\ \ff t\in[0,T].\end{equation}By \eqref{e3.10}, \eqref{e3.11}, it follows that there is an adapted process
\begin{equation}\label{e3.11a}
\barr{l}
y\in M^2_\PP(0,T;H^1_0(\calo)\cap L^2_{ad}(\ooo;C([0,T];L^2(\calo))\vsp
\dd\frac{\pp y}{\pp t}\in M^2_\PP(0,T);L^2(\calo)),\earr\end{equation} such that, for $\vp\to0$,
$$\barr{lcll}
y_\vp&\to&y& \mbox{strongly in $L^2_{ad}(\ooo;C([0,T];L^2(\calo))\cap M^2_\PP(0,T;H^1_0(\calo)),$}\earr$$
and $\pas$,
\begin{equation}\label{e3.12}
\barr{rcll}
\D y_\vp&\to&\D y& \mbox{weakly in $L^2(Q_T),\ \pas,$}\vsp
\b_\vp(y_\vp)&\to&\eta& \mbox{weakly in $L^2(Q_T),\ \pas,$}\vsp
\wt F(t,y_\vp)&\to&\wt F(t,y)&\mbox{strongly in $L^2(Q_T),$}\vsp
\dd\frac{\pp y_\vp}{\pp t}&\to&\dd\frac{\pp y}{\pp t}&\mbox{weakly in $L^2(Q_T)$},\earr\end{equation}where $\eta\in\b(y)$, a.e. in $Q$, $\pas$ We have also $\pas$
\begin{equation}\label{e3.13}
\barr{l}
\dd\frac{\pp y}{\pp t}-\D y+\wt F(t,y)+g\cdot\na y+\eta=f,\mbox{\ \ a.e. in }Q_T,\vsp
y(0,\xi)=x(\xi),\ \inr\ \calo,\qq y=0\ \onr\ \Sigma_T.\earr\end{equation}Since $y$ is a smooth $y\in W^{1,2}([0,T];L^2(\calo))$ and predictible solution to \eqref{e3.13}, we conclude that $X$ given by \eqref{e3.1} is the strong solution to \eqref{e2.1} (see \cite{9}, p.~217).  The uniqueness of the solution $X$ follows by a standard argument involving Ito's formula. Moreover, since $$\E\{\sup(e^\mu+|\D\mu|;\ t\in[0,T])^n\}<\9,\ \ff n\in\nn,$$ we infer by \eqref{e3.11a} that \eqref{e2.7a} holds. This completes the proof.

\subsection{Proof of Theorem \ref{t2.4}}

Proceeding as in the proof of Theorem \ref{t2.2}, we reduce, via transformation \eqref{e3.1}, equation \eqref{e1.5} to the random equations
\begin{equation}\label{e3.14}
\left\{\barr{ll}
\dd\frac{\pp y}{\pp t}-\D y+\wt F(t,y)+g\cdot\na y=f&\inr\ Q_T,\vsp
y(0)=x&\inr\ \calo,\vsp
\dd\frac{\pp y}{\pp\nu}+\frac{\pp\mu}{\pp\nu}\, y+\b(y)\ni0&\onr\ \Sigma_T,\earr\right.\end{equation}where $\wt F$ and $g$ are given by \eqref{e3.4}. Now, we fix $\oo\in\ooo$ and treat \eqref{e3.14} as a deterministic nonlinear parabolic \eq.

It must be said that, though the existence theory for nonlinear parabolic boundary value problems is not directly applicable in the present situation, we may prove, however, existence in \eqref{e3.14} via standard approximation technique. Namely, we approximate \eqref{e3.14} by the penalized equation
\begin{equation}\label{e3.15}
\left\{\barr{ll}
\dd\frac{\pp y_\vp}{\pp t}-\D y_\vp+\wt F(t,y_\vp)+g\cdot\na y_\vp=f&\inr\ Q_T,\vsp
y_\vp(0)=x&\inr\ \calo,\vsp
\dd\frac{\pp y_\vp}{\pp \nu}+\dd\frac{\pp \mu}{\pp \nu}\,y_\vp+\b_\vp(y_\vp)=0&\onr\ \Sigma_T.\earr\right.\end{equation}Equivalently,
\begin{equation}\label{e3.16}
\barr{l}
\dd\frac{dy_\vp}{dt}+\wt A_\vp(t)y_\vp=f(t)\q\inr\ (0,T),\vsp
y_\vp(0)=x,\earr\end{equation}where $\wt A_\vp(t):V\to V'$ is defined by
\begin{equation}\label{e3.17}\barr{lcl}
_{V'}\<\wt A_\vp(t)\vf,\vf\>_V&=&\dd\int_\calo(\na y_\vp\cdot\na\vf+\wt F(t,y_\vp)\vf+(g\cdot\na y_\vp)\vf)d\xi\vsp
&&+\dd\int_{\pp\calo}\(\b_\vp(y_\vp) +\dd\frac{\pp\mu}{\pp\nu}\,y_\vp \)\vf\,d\xi,\ \ff\vf\in V.\earr\end{equation}We have, for all $t\in[0,T]$,
\begin{eqnarray}
\|\wt A_\vp(t)y\|_{V'}&\le&C_1\|y\|_V,\q\ff y\in V,\q\label{e3.18}\\[1mm]
_{V'}\<\wt A_\vp(t)y,y\>_V&\ge&C_2\|y\|^2_V-C_3\|y\|^2_H,\q\ff y\in V,\q\label{e3.19}\\[1mm]
_{V'}(\wt A_\vp(t)y-\wt A_\vp(t)\bar y,y-\bar y)_V&\ge&-C_4\|y-\bar y\|^2_H,\q\ff y\in V,\q\label{e3.20}
\end{eqnarray}
where $C_i$, $i=1,..,4,$ are independent of $\vp$.
To get \eqref{e3.19}, \eqref{e3.20}, we have used the trace-interpolation estimate
$$\left|\int_{\pp\calo}\dd\frac{\pp\mu}{\pp\nu}\,y^2_\vp\,d\xi\right|\le C\|y_\vp\|^2_{H^{\frac12}(\pp\calo)}\le C\|y_\vp\|_{L^2(\calo)}\|y_\vp\|_{H^1(\calo)}.$$

Then, by a standard existence result for the Cauchy problem \eqref{e3.16} (see J.L. Lions \cite{12} or \cite{1}, p.~177), there is a unique $y_\vp\in L^2(0,T;V)\cap C([0,T];H)$, with $\dd\frac{y_\vp}{dt}\in L^2(0,T;V')$ which satisfies \eqref{e3.16}, a.e. on $(0,T)$. It should be mentioned that, since the solution $y_\vp$ can be obtained via Galerkin scheme or by iteration via Banach contraction principle, the process $y_\vp=y_\vp(t,\oo)$ is adapted in the stochastic probability basis $\{\ooo,\PP,\calf,\{\calf_t\}_{t\ge0}\}$.

Now, by \eqref{e3.16} or by \eqref{e3.16}, \eqref{e3.17}, we have the \fwg\ a priori estimates
$$\barr{r}
\dd\frac12\ \dd\frac d{dt}\,\|y_\vp(t)\|^2_H+\|\na y_\vp(t)\|^2_H\le C(\|y_\vp(t)\|^2_H+\|y_\vp(t)\|^2_{L^2(\pp\calo)})+\|f(t)\|^2_H,
\vsp\mbox{ a.e. }t\in(0,T).\earr$$This yields, as in the proof of Theorem \ref{t2.2},
\begin{equation}\label{e3.17a}
\E\|y_\vp(t)\|^2_H+\E\int^t_0\|y_\vp(s)\|^2_Vds\le C,\q t\in[0,T].\end{equation}(We shall denote by $C$ several positive constants independent of $\vp.$)

Now, we take in \eqref{e3.3} $\vf=\b_\vp(y_\vp)$ and obtain that
$$\barr{r}\dd\frac d{dt}\int_\calo j_\vp(y_\vp(t,\xi))d\xi+\dd\int_{\pp\calo}|\b_\vp(y_\vp(t,\xi))|^2d\xi\le C(\|y_\vp(t)\|^2_V+\|f(t)\|^2_H),\vsp\mbox{ a.e. }t\in(0,T),\earr$$where $j_\vp(r)=\dd\int^r_0\b_\vp(s)ds.$ This yields
\begin{equation}\label{e3.18a}
\barr{l}
\dd\int_\calo j_\vp(y_\vp(t,\xi))d\xi+\dd\int^t_0\int_{\pp\calo}
\b^2_\vp(y_\vp(s,\xi))d\xi\vsp\qq\le \dd\int_\calo j_\vp(x)d\xi+C\dd\int^t_0\|f(s)\|^2_Hds\le C,\ \ff \vp>0,\ t\in[0,T].\earr\end{equation}Finally, if in \eqref{e3.17a} we take $\vf=y_\vp-y_{\vp'}$ and use \eqref{e3.16}, we obtain that
$$\barr{l}
\dd\frac12\,\|y_\vp(t)-y_{\vp'}(t)\|^2+\dd\frac12\int_0^t\|y_\vp(s)-y_{\vp'}(s)\|^2_Vds\vsp
\qq+\dd\int^t_0\int_{\pp\calo}(\b_\vp(y_\vp(s,\xi))-\b_{\vp'}(y_{\vp'}(s,\xi)))(y_\vp(s,\xi)-y_{\vp'}(s,\xi))ds\,d\xi\vsp
\qq\le C(\oo)\dd\int^t_0\|y_\vp(s)-y_{\vp'}(s)\|^2_Hds,\q \ff t\in[0,T].\earr$$
Then, by \eqref{e3.9} we see that
$$\|y_\vp(t)-y_{\vp'}(t)\|^2_H+\int^t_0\|y_\vp(s)-y_{\vp'}(s)\|^2_Vds\le C(\vp+\vp')y,\ \pas$$and, therefore,
\begin{equation}\label{e3.19a}
\barr{lcll}y_\vp&\to&y&\inr\ L^2(0,T;V)\cap C([0,T];H)\earr\end{equation}as $\vp\to0$. Note, also, that the process $y=y(t,\oo)$ is adapted. Moreover, since, by \eqref{e3.9}, $\{\b_\vp(y_\vp)\}$ is bounded in $L^2(\Sigma)$, we have (eventually on a subsequence)
\begin{equation}\label{e3.20a}
\barr{lcll}\b_\vp(y_\vp)&\to&\eta&\mbox{weakly in } L^2(\Sigma_T).\earr\end{equation}On the other hand, by \eqref{e3.19a} and the trace theorem,  it follows that, for $\vp\to0$,
\begin{equation}\label{e3.21a}
\barr{lcll}y_\vp&\to&y&\mbox{strongly in } L^2(\Sigma_T)\earr.\end{equation}Since $\b$ is maximal monotone in $L^2(\Sigma_T)\times L^2(\Sigma_T)$, we infer by \eqref{e3.20a} and \eqref{e3.21a} that
\begin{equation}\label{e3.22a}
\eta(t,\xi)\in\b(y(t,\xi)),\mbox{\ \ a.e. }(t,\xi)\in\Sigma_T.\end{equation}Then, letting $\vp$ tend to zero into \eqref{e3.3}, we see that
\begin{equation}\label{e3.23a}
\barr{l}
\dd\frac{dy}{dt}+\wt A(t)y=f(t),\q\mbox{ a.e. }t\in(0,T),\vsp
y(0)=x,\earr\end{equation}where $\wt A(t)$ is given by
\begin{equation}\label{e3.24a}\barr{lcl}
_{V'}(\wt A(t)y,\vf)_V&=&\dd\int_\calo(\na_\xi y)\cdot\na\vf+\wt F(t,y)\vf+(g\cdot\na y)\vf)d\xi\vsp&&+\dd\int_{\pp\calo}\(\eta+\dd\frac{\pp\mu}{\pp\nu}\,y\)\vf\,d\xi,\ \ff\vf\in V.\earr\end{equation}

Now, if we take $X=e^\mu y$ as in \eqref{e3.1}, taking into account that $y$ is a regular predictible solution to \eqref{e3.23a}, we conclude   that $X$ is   a strong solution to \eqref{e1.6} in the sense of Definition \ref{d2.3}. By \eqref{e3.18a} it follows also that $$y(t,\xi,\oo)\ge0,\mbox{\ \ a.e. }(t,\xi,\oo)\in Q_T\times\ooo,$$and, therefore, $X\ge0$, a.e. in $Q_T\times\ooo$. The uniqueness follows by \eqref{e2.6} (or~\eqref{e2.7}) via Ito's formula, but the details are omitted. This completes the proof.

\section{The stochastic one-phase Stefan problem}
\setcounter{equation}{0}

The Stefan problem describes the heat conduction in a medium $\calo\subset\rr^d$, $d=1,2,3$, phase change. Assume that, at time $t\ge0$, $\calo^t\subset\calo$ is the solid region (ice) and that $\calo^t=\{\xi\in\calo;\ 0<\ell(\xi)<t<T\},$ $\pp\calo^t=\{t=\ell(\xi)\}$, where $\ell$ is a smooth function. Then the classical one-phase Stefan problem is described by the system (see, e.g., \cite{1}, p.~231)
\begin{equation}\label{e4.1}
\left\{\barr{ll}
\dd\frac{\pp\theta}{\pp t}-\D\theta=0&\inr\ \{(t,\xi)\in Q_T;\ \ell(\xi)<t<T\},\vsp
\theta=0&\inr\ \{(t,\xi);\ t<\ell(\xi)\},\vsp
\na_\xi\theta\cdot\na\ell=-\rho&\onr\ S=\{(t,\xi);\ t=\ell(\xi)\},\vsp
\theta(0,\xi)=\theta_0(\xi)\ \inr\ \calo,&\theta=0\ \onr\ \Sigma_T,\earr\right.
\end{equation}where $\theta=\theta(t,\xi)$ is the temperature. Under a stochastic multiplicative perturbation of the form $W=\theta\sumk\mu_k\b_k$, system \eqref{e4.1} becomes
\begin{equation}\label{e4.2}
\barr{ll}
d\theta-\D\theta\,dt=\theta\sumk\mu_kd\b_k&\inr\ \{(t,\xi)\in Q_T;\ \ell(\xi)<t<T\},\vsp
\theta=0&\inr\ \{(t,\xi);\ 0<t<\ell(\xi)\},\vsp
\na_\xi\theta\cdot\na\ell=-\rho&\onr\ S,\vsp
\theta(0)=\theta_0\ \inr\ \calo,&\theta=0\ \onr\ \Sigma.\earr\end{equation}Roughly speaking, this means that in  (melting) processes arises a heating noise source of the form $\sumk\theta\mu_k\dot\b_k$.

We shall show below that \eq\ \eqref{e4.2} reduces to an obstacle stochastic problem of the form \eqref{e1.1}. To this end, it is convenient to reduce \eqref{e4.2} to a random one-phase Stefan problem via the transformation
\begin{equation}\label{e4.3}
\theta=e^\mu z\q\inr\ Q_T.\end{equation}Arguing as in the proof of Theorem \ref{t2.2}, we see that this happens if $z$ is solution to the equation
\begin{equation}\label{e4.4}
\left\{\barr{l}
\dd\frac{\pp z}{\pp t}-\D z-2\na\mu\cdot\na z+(\wt\mu-\D\mu-|\na\mu|^2)z=0\vsp
\hfill\inr\ \{(t,\xi)\in Q_T;\ \ell(\xi)<t<T\},\vsp
z=0\ \inr\ \{(t,\xi);\ 0<t\le\ell(\xi)\},\vsp
\na_\xi  z\cdot\na\ell=-\rho\ \onr\ S,\vsp
z(0)=\theta_0\ \inr\ \calo,\qq z=0\ \onr\ \Sigma_T.\earr\right.\end{equation}Equation \eqref{e4.4} is a random one-phase free boundary parabolic equation of Stefan type which can be reduced to a parabolic variational inequality via the transformation
$$y(t,\xi)=\left\{\barr{ll}
\dd\int^t_{\ell(\xi)}\theta(s,\xi)ds&\ifr\ \xi\in\calo\setminus\calo^0,\ t>\ell(\xi),\vsp
\dd\int^t_0\theta(s,\xi)ds&\ifr\ \xi\in\calo^0,\ 0<t<T.\earr\right.$$We set
$$f(t,\xi)=\left\{\barr{ll}
-\rho&\ifr\ \xi\in \calo\setminus\calo^0,\ 0<t<T,\vsp
\theta_0(\xi)&\ifr\ \xi\in\calo^0,\ 0<t<T,\earr\right.$$where $\calo^0=\{\xi\in\calo;\ \theta_0(\xi)>0\}.$

Then, arguing as in Lemma 5.1 in \cite{1}, it follows that $y$ is the solution to the \eq
\begin{equation}\label{e4.5}
\barr{l}
\dd\frac{\pp y}{\pp t}-\D y-2\na\mu\cdot\na y+(\wt\mu-\D\mu-|\na\mu|^2)y=f_0\dubluunu_{Q_+}\ \inr\ \cald'(Q),\vsp
y(0)=0\ \inr\ \calo;\qq y=0\ \onr\ \Sigma,\earr\end{equation}
where $Q_+=\{(t,\xi);\ \ell(\xi)<t<T\}$ and $\cald'(Q)$ is the space of distributions on $Q$. 
Equivalently,
\begin{equation}\label{e4.6}
\barr{l}
\dd\frac{\pp y}{\pp t}-\D y-2\na\mu\cdot\na y+(\wt\mu-\D\mu-|\na\mu|^2)y+\b(y)\ni f_0\dubluunu_{Q_+}\vsp
y(0)=0\ \inr\ \calo;\qq y=0\ \onr\ \Sigma_T,\earr\end{equation}
and so, $X=e^\mu y$ is solution to equation \eqref{e1.1}, where $f=f_0\dubluunu_{Q_+}e^\mu$. In this way, the existence into problem \eqref{e4.2} reduces to Theorem \ref{t2.2}.

\end{document}